\documentclass[10pt]{article}
\usepackage{latexsym}
\usepackage{amssymb}
\usepackage{amsmath}
\usepackage{layout}
\newtheorem{teo}{Theorem}[section]
\newtheorem{lema}[teo]{Lemma}

\newtheorem{prop}[teo]{Proposition}
\newtheorem{obs}[teo]{Remark}
\newtheorem{defin}{Definition}
\newenvironment{dem}{\text\it Proof:}{}
\newcommand{\R}{\mathbb{R}}
\newcommand{\N}{\mathbb{N}}

\newcommand{\av}{\mbox{{\bf Av}}}
\newcommand{\be}{\mbox{{\bf E}}}
\newcommand{\bm}{\beta_{-}}
\newcommand{\ent}{{\cal E}_{n,t}}

\newcommand{\ou}{[0,1]}
\newcommand{\nut}{\nu_t}
\newcommand{\sii}{\sigma_i}

\newcommand{\ssm}{\Sigma_{N-1}}
\newcommand{\ssn}{\Sigma_N}

\newcommand{\ce}{\mathcal E}
\newcommand{\cf}{\mathcal F}

\newcommand{\ep}{\varepsilon}

\newcommand{\si}{\sigma}

\newcommand{\vp}{\varphi}

\newcommand{\lp}{\left(}
\newcommand{\rp}{\right)}

\newcommand{\lla}{\left\langle}
\newcommand{\rra}{\right\rangle}

\newcommand{\ram}{\right\rangle_{-}}

\newcommand{\rat}{\right\rangle_{t}}

\begin{document}

 \thispagestyle{empty}
 \begin{center}
 {\Large\bf Higher order expansions}\\
{\Large \bf for the overlap
of the SK model}

\vspace{0.3cm}

 by

\vspace{0.3cm}

{\bf Xavier Bardina}\footnote{ Partially supported
by DGES grants BFM2000-0009, BFM2000-0607,  HF2000-0002.}, {\bf David M\'arquez-Carreras$^2$, Carles Rovira}\footnote{{Partially supported by DGES grants
BFM2000-0607 and HF2000-0002
.\hfill}}  and {\bf Samy
Tindel}$^3$

\vspace{0.5cm}

$^1$
{\it Departament de Matem\`atiques, Universitat Aut\`onoma de Barcelona}
\\
\it 08193 Bellaterra, Barcelona, Spain
\\
{\it e-mail: bardina@mat.uab.es}
\\
 $^2$
{\it Facultat de Matem\`atiques,
Universitat de Barcelona,}
\\
\it Gran Via 585,
08007-Barcelona,
Spain
\\
{\it e-mail: marquez@mat.ub.es, rovira@mat.ub.es}
\\
$^3$ {\it D\'epartement de Math\'ematiques,
Institut Galil\'ee - Universit\'e Paris 13,}
\\
{\it Avenue  J. B. Cl\'ement,
93430-Villetaneuse,
France}
\\
{\it e-mail: tindel@math.univ-paris13.fr}

 \end{center}

\vspace{0.2cm}

\begin{abstract}In this note, the Sherrington Kirkpatrick model of interacting spins is under consideration. In the high temperature region, we give an asymptotic expansion for the expected value of some genereral polynomial of the overlap of the system when the size $N$ grows to infinity. Some of the coefficients obtained are shown to be vanishing, while the procedure to get the nontrivial ones has to be performed by a computer program, due to the great amount of computation involved.

\vspace{0.2cm}

\end{abstract}

\section{Introduction}
This note is concerned with the usual Sherrington-Kirkpatrick model of interacting spins, that can be briefly described as follows: for $N\ge 1$, set $\ssn=\{ -1,1 \}^N$, called the space of configurations. For a given $\si=(\si_1,\ldots,\si_N)\in\ssn$, the quantity $\sii$ will stand for the value of the $i^{\mbox{\tiny th}}$ spin of the configuration. We will suppose that the spins randomly interact according to a mean field model represented by a Hamiltonian $H_N$ given by
$$
-H_N(\si)=\frac{1}{N^{1/2}}\sum_{1\le i<j\le N}g_{i,j}\sii \si_j,
$$
where $\{g_{i,j}; 1\le i<j \le N  \}$ is a family of independent standard Gaussian random variables, defined on a probability space $(\Omega,\cf,{\bf P})$. For a given positive parameter $\beta$ (that stands for the inverse of the temperature of the system), set then
$$
Z_N=\sum_{\si\in\ssn}\exp\lp-\beta H_N(\si)  \rp.
$$
We will consider here the random Gibbs measure $G_n$ on $\ssn$, that gives the highest weight to those configurations minimizing the energy $H_N$, and whose density with respect to the uniform measure is given by $Z_N^{-1}\exp(-\beta H_N(\si))$ for all $\si\in\ssn$. For any $N,n\ge 1$ and any function $f:\ssn^n\to\R$, we will then set $\langle f\rangle$ for the average of $f$ with respect to $dG_N^{\otimes n}$, i.e.
$$
\lla f \rra =Z_N^{-n}
\sum_{\si^1,\ldots,\si^n\in\ssn}f(\si^1,\ldots,\si^n)
\exp\lp-\beta \sum_{l\le n}H_N(\si^l)  \rp,
$$
and the quantity $\be\langle f\rangle$ will be denoted by $\nu(f)$. 

Parisi's analysis of the SK model (see \cite{MPV} for a complete discussion) showed that
a family of functions of particular interest, summarizing a great amount of information about the whole spin system, is the family $\{R_{l,l'};1\le l<l'<\infty  \}$, where
$$
R_{l,l'}=\frac{1}{N}\sum_{i\le N}\sii^l \sii^{l'},
$$
and $\si^l, \si^{l'}$ are understood as two independent configurations under the measure $dG_N^{\otimes 2}$. The quantity $R_{l,l'}$ (which will be the main object of study of the present paper) is called overlap between the configurations $\si^l$ and $\si^{l'}$, and can be easily related to the number of spins of $\si^l$ and $\si^{l'}$ which are in the same direction.

While the efforts of Physicits generally aim at the understanding of the low temperature regime of the SK model, and of the related phenomenon of non self averaging (see e.g. \cite{PS} and \cite{T98}), the high temperature domain (i.e. $\beta\in[0,1)$) is still of mathematical interest: indeed, in that case, a rather complete description of the limiting behavior of the system when  $N\to\infty$ can be given, generally amounting to some asymptotic theorems concerning $Z_N$ and $R_{l,l'}$. Among those results, we can cite \cite{Tsf} concerning the $L^2$-limit of $\frac{1}{N}\log(Z_N)$, as well as \cite{AZL}, \cite{CN} and \cite{GT}, that provide some Central Limit Theorems for $Z_N$ and $R_{1,2}$. However, the possibilities of getting some sharp estimates for the above quantities have been greatly enhanced by the rigorous introduction and proof of the cavity method, summarized in \cite{Tbk}, leading among other results to a number of interesting CLT for the family $\{R_{l,l'};1\le l<l'<\infty  \}$.

In this paper, we would like to go one step further into that line of investigation, and try to give a complete expansion (in terms of power series of $\frac{1}{N}$) for the quantity $\nu(R_{1,2}^2)$, and more generally, for any quantity of the form
\begin{equation}\label{zobi}
\nu\lp\prod_{i=1}^{r}\si_N^{j_i}\si_N^{l_i}R_{j_i,l_i}  \rp,
\end{equation}
for an arbitrary sequence of indexes $(j_i,l_i)$. This is indeed a natural question, quoted as Research Problem 2.3.1 in \cite{Tbk}. Furthermore,  it can be shown, using the stochastic calculus methods introduced in \cite{CN}, as well as the large deviations estimates for semi-martingales given by \cite{Pu}, that a moderate deviation principle could be obtained for $\frac{1}{N}\log(Z_N)$ as long as we can prove the boundedness in $N$ of some exponential moments for the random variable
$$
N\lp\lla NR_{1,2}^2 \rra-\frac{1}{1-\beta^2}\rp.
$$
This, in turn, can be easily reduced to the asymptotic behavior of quantitites of the form
$$
\lla \prod_{j=1}^{k}R^2_{2j-1,2j} \rra,
$$
for any $k\ge 1$. This article can be understood then as a preliminary step towards the moderate deviations principle for $Z_N$.

Let us now present the main results of our paper. The general expansion in terms of $N$ for (\ref{zobi}) is given by the following theorem:
\begin{teo}\label{teo1}
For any $n \ge \left[\frac{k+1}2\right]$ and $\beta<1$, there exist constants $C(\beta,j)$ such that
$$\nu(\prod_{i=1}^k \si_N^{j_i}\si_N^{l_i}R_{j_i,l_i})=\sum_{j=\left[\frac{k+1}2\right]}^n\frac{C(\beta,j)}{N^j}
+O\left(\frac1{N^{n+1}}\right).$$
\end{teo}

\begin{obs}
One would also expect to get some terms of order $N^{\frac{2j+1}{2}}$ in that expansion. Our result shows that the coefficients corresponding to those terms are all vanishing.
\end{obs}

As a particular case of our general result, we have the following

\begin{teo}\label{teo2}
For any $n \ge 1$ and $\beta<1$, there exist constants $C(\beta,j)$ such that
$$\nu(R_{1,2}^2)=\sum_{j=1}^n\frac{C(\beta,j)}{N^j}
+O\left(\frac1{N^{n+1}}\right),$$
where
$$
C(\beta,1)=\frac{1}{(1-\beta^2)},
C(\beta,2)=\frac{-\beta^2(1+\beta^2)}{(1-\beta^2)^4},
C(\beta,3)=\frac{\frac{2}{3} (6+19\beta^2+4\beta^4-5\beta^6)}{(1-\beta^2)^7}.
$$
\end{teo}

In fact, Theorem \ref{teo1} will be an easy consequence of the next proposition, whose proof will be detailed in the sequel of the paper.
\begin{prop}\label{propgen}
For any $n \ge  \left[\frac{k+1}2\right]$ and $\beta<1$, a term ${\displaystyle
\nu(\prod_{i=1}^k\si_N^{j_i}\si_N^{l_i}R_{j_i,l_i})}$
can be written as a finite sum of terms of the type
$$K(\beta,I)\frac1{N^j}\nu(\prod_{i=1}^r\si_N^{j_i}\si_N^{l_i}R_{j_i,l_i})$$
with $j,r \ge 0$, $j+\frac r2\geq \frac k2$ plus a remainder of order
$O\left(\frac1{N^{n+\frac{1}{2} }}\right)$, where $K(\beta,I)$ denotes a constant depending on $\beta$ and
the family of indexes $I=\{(j_i,l_i), 1 \le i \le r \}$. Moreover, when $j+\frac r2=\frac k2$ then $j>0$.
\end{prop}

\begin{obs}
Notice that when $r=0$ we get the terms of the form $\frac{K(\beta,I)}{N^j}$ with $j \in \N$
and $j\geq \frac k2$. These terms are the unique ones that appear in the decomposition given
in Theorem \ref{teo1}.
\end{obs}

Note that the bulk of the computations one has to perform in order to obtain the expansion of $\nu(R_{1,2}^2)$ does not allow us to give a closed formula for all the constants $C(\beta,j)$. However, the analysis of the mechanism needed to obtain those fomulae lead us to write a computer program performing this kind of calculations. The listing of the program will be given in an Appendix.

Our paper will be divided as follows: in the next Section, we will recall some basic facts and methods concerning the SK model. At Section 3, for sake of readability, we will explain our methods by calculating the first two terms of the expression of $\nu(R_{1,2}^2)$. Then, at Section 4, we will prove our main result (i.e. Proposition \ref{propgen}). Section 5 will be devoted to the explanation of the algorithm of our program, which will be given in the Appendix.

Troughout the paper all constants appearing in the proofs are called $K, K(\beta), K(\beta,k)$ and $K(\beta,k,n)$,
although they may vary from one part to another. We will omitt the dependence on the family of indexes $I$.

\section{The cavity method}
Let us recall now the basic features of the cavity method (taken from \cite{Tbk}), that allows to reduce integrals with respect to a $N$-spins system to a $N-1$-spins one: set
\begin{equation}\label{defbetaprime}
\bm=\lp \frac{N-1}{N} \rp^{\frac{1}{2}}\beta,
\end{equation}
and $\langle . \rangle_{-}$ the averaging with respect to $G_N$ at temperature $\bm$. Let $n\ge 1$ and $\si^1,\ldots,\si^n$ be $n$ independent copies of a $N$-spins configuration. For any $j \in \{1,\ldots,n\},$ we will also denote $\si^j=(\rho^j,\ep^j)$, where $\rho^j\in\ssm$ and $\ep^j\in\{-1,1  \}$. Set
\begin{eqnarray}
g_i&=&g_{i,N}\nonumber\\
g\lp \rho\rp&=&\frac{\beta}{N^{1/2}}\sum_{i\le N-1}g_i\sii.\nonumber
\end{eqnarray}
Then it is easily shown that for one configuration $\si$, setting $\si_N=\ep$, we have
$$
-\beta H_{N}(\si)=-\bm H_{N-1}(\rho)+\ep g\lp \rho\rp,
$$
and hence, for $f:\ssn^n\to\R$, we can decompose $\langle f \rangle$ into
$$
\lla f \rra= \frac{\lla \av f \ce_n  \ram}{Z^{n}},
$$
where, for a function $\vp:\{-1,1  \}^n\to\R$, we set
$$
\av \vp =\frac{1}{2^n}\sum_{\ep_1,\ldots,\ep_n\in\{ -1,1 \}^n}\vp(\ep_1,\ldots,\ep_n),
$$
and
\begin{eqnarray*}
\ce_n&=&\exp\lp \sum_{l\le n}\ep^l  g(\rho^l) \rp\\
Z&=&\lla\av \ce_1\ram=\lla\cosh\lp g(\rho^l) \rp\ram.
\end{eqnarray*}

The path between the original spins system and a system in which the last spin is decoupled from the other ones, can be defined as follows: for $t\in\ou$, set
$$
g_t(\rho)=t^{1/2}g(\rho).
$$
For $n\ge 1$, $f:\ssn^n\to\R$ and $t\in\ou$, let then $\langle f \rangle_t$ be the quantity defined by
$$
\lla f\rat=Z_t^{-n}\lla \av f \ent \ram,
$$
where
\begin{eqnarray*}
\ent&=&\exp\lp \sum_{l\le n}\ep^l  g_t(\rho^l) \rp\\
Z_t&=&\lla \av \ce_{1,t}\ram=\lla\cosh\lp g_t(\rho^l)) \rp\ram.
\end{eqnarray*}
Set also $\nut(f)=E[\langle f \rangle_t]$ for $t\in\ou$, and notice that $\nu(f)=\nu_1(f)$.

The following propositions, that we label for further use, are then shown in \cite{Tbk}.
\begin{prop}\label{deriv}
For a function $f$ on $\Sigma_N^n$, and $t>0$, we have
\begin{eqnarray*}
\nu'_t(f)&=&\beta^2 \sum_{1 \le l < l' \le n} \nu_t(f \varepsilon_l \varepsilon_{l'} R_{l,l'})
\\
&-&\beta^2 n\sum_{l \le n} \nu_t(f \varepsilon_l \varepsilon_{n+1} R_{l,n+1})
\\
&+&\beta^2 \frac{n(n+1)}{2} \nu_t(f \varepsilon_{n+1} \varepsilon_{n+2} R_{n+1,n+2}).
\end{eqnarray*}
\end{prop}

\begin{prop}\label{der1}
If $f$ is a function on $\Sigma_N^n$, then
$$\vert \nu_t^{(k)}(f)\vert \le \frac{K(\beta,k,n)}{N^\frac{k}{2}} \nu(f^2)^{\frac12}.$$
\end{prop}

\section{Expansion of $\nu(R_{1,2}^2)$}

We begin by introducing some notation.

\begin{defin} Given an expression $\prod_{i=1}^k\varepsilon_{j_i}\varepsilon_{l_i}R_{j_i,l_i}$  we define the
$\varepsilon$-order of this expression as the number of different indexes that appears an odd number of times in the set $\{j_i,l_i;\, 1 \le i \le k \}$. Notice that the   $\varepsilon$-order can be computed as
$$ \sum_n \Big( \mod \big( \sum_{i=1}^k (1_{\{l_i=n\}} + 1_{\{j_i =n\}} ),2 \big) \Big).$$
\end{defin}

Notice that the $\varepsilon$-order is always a nonnegative even number. Note also that
$$
\prod_{i=1}^k\varepsilon_{j_i}\varepsilon_{l_i}R_{j_i,l_i}
\quad\mbox{ and } \quad
\prod_{i=1}^k\varepsilon_{j_i}\varepsilon_{l_i}R_{j_i,l_i}^{-}
$$
are of the same $\varepsilon$-order.

\medskip
In order to get the expansion of $\nu(R_{1,2}^2)$ we need the following technical lemma that
allows us to obtain of a lot of terms.

\begin{lema}\label{simp1}
Given ${\displaystyle\prod_{i=1}^k
\varepsilon_{j_i}\varepsilon_{l_i}R_{j_i,l_i}}$, set $\lambda$ its corresponding $\varepsilon$-order. Then for all $\beta <1$
\begin{eqnarray*}
&(a)&
\nu_0(\prod_{i=1}^k\varepsilon_{j_i}\varepsilon_{l_i}R_{j_i,l_i}^{-})=0
\quad
\textrm{when $\lambda>0$, }\\
&&|\nu_0(\prod_{i=1}^k\varepsilon_{j_i}\varepsilon_{l_i}R_{j_i,l_i}^{-})|\leq\frac{K(k,\beta)}{N^{\frac
k2}} \quad \textrm{when $\lambda=0$,}\\
&(b)&
|\nu(\prod_{i=1}^k\varepsilon_{j_i}\varepsilon_{l_i}R_{j_i,l_i}^{-})|\leq\frac{K(k,\beta)}{N^{\frac
k2 },}\\
&(c)&|\nu(\prod_{i=1}^k\varepsilon_{j_i}\varepsilon_{l_i}R_{j_i,l_i})|\leq\frac{K(k,\beta)}{N^{\frac
k2 }.}
\end{eqnarray*}
\end{lema}
\begin{dem}
The case (a) when $\lambda>0$ follows immediately from the definition of $\nu_0$.

The other cases follow easily from the estimates $\nu(R_{1,2}^{2k}) \le \frac{K(k,\beta)}{N^k}$,
 $\nu((R^-_{1,2})^{2k}) \le \frac{K(k,\beta)}{N^k}$ and
$\nu_t((R^-_{1,2})^{2k}) \le K(k,\beta) \nu((R^-_{1,2})^{2k}) $ given by Talagrand (see (2.88), (2.89) and (2.91) in \cite{Tbk}), and the following version of H\"older's inequality: if $\tau_1, \tau_2>1, \frac{1}{\tau_1} + \frac{1}{\tau_2}=1$, then for all $t \in [0,1]$
$$\nu_t(f_1 f_2) \le \nu_t( \vert f_1 \vert^{\tau_1})^{\frac{1}{\tau_1}}
\nu_t( \vert f_2 \vert^{\tau_2})^{\frac{1}{\tau_2}}.$$
\hfill $\Box$
\end{dem}
\bigskip

Then we are able to prove the following expansion.

\begin{prop}
For all $\beta <1$ we have
$$
\nu(R_{1,2}^2)=\frac1{N(1-\beta^2)}-\frac{\beta^2(1+\beta^2)}{N^2(1-\beta^2)^4}
+O\left(\frac1{N^{\frac52}}\right).
$$
\end{prop}

\begin{dem}
Using the symmetry property and the relation $\varepsilon_1 \varepsilon_2 R_{1,2} =\frac1N + \varepsilon_1\varepsilon_2R_{1,2}^-$ we get
\begin{equation}\label{p1}
\nu(R_{1,2}^2)=\nu(\varepsilon_1\varepsilon_2R_{1,2})=\frac1N+\nu(\varepsilon_1\varepsilon_2R_{1,2}^{-}),
\end{equation}
Performing a Taylor expansion and using Proposition \ref{der1} and Lemma \ref{simp1} in order to bound the last term we get
\begin{eqnarray}
\nu(\varepsilon_1\varepsilon_2R_{1,2}^{-})&=&\nu_0(\varepsilon_1\varepsilon_2R_{1,2}^{-})+
{\nu'\!}_0(\varepsilon_1\varepsilon_2R_{1,2}^{-})\nonumber\\
&&+\frac12{\nu''\!}_0(\varepsilon_1\varepsilon_2R_{1,2}^{-})+\frac16
\nu_0^{(3)}(\varepsilon_1\varepsilon_2R_{1,2}^{-})+O\left(\frac1{N^{\frac52}}\right)\label{p2}.
\end{eqnarray}
Notice first that $\nu_0(\varepsilon_1\varepsilon_2R_{1,2}^{-})=0$. To deal with the $k$-derivative term we will apply Proposition \ref{deriv} $k$ times iteratively and we simplify all the null terms using the first part of Lemma \ref{simp1}.

{\it Study of ${\nu'\!}_0(\varepsilon_1\varepsilon_2R_{1,2}^{-})$}. It is easy to check that
$${\nu'\!}_0(\varepsilon_1\varepsilon_2R_{1,2}^{-})=\beta^2\nu_0((R_{1,2}^{-})^2).$$
Using again a Taylor's expansion we can write
\begin{eqnarray}
{\nu'\!}_0(\varepsilon_1\varepsilon_2R_{1,2}^{-})&=&\beta^2\nu_0((R_{1,2}^{-})^2)\nonumber\\
&=&\beta^2\left[\nu((R_{1,2}^{-})^2)-{\nu'\!}_0((R_{1,2}^{-})^2)-\frac12{\nu''\!}_0((R_{1,2}^{-})^2)+O\left(\frac1{N^{\frac52}}\right) \right].\label{p3}
\end{eqnarray}
From Proposition \ref{deriv} and Lemma \ref{simp1} we get that ${\nu'\!}_0((R_{1,2}^{-})^2)=0$. On the other hand, using again the relation $\varepsilon_1 \varepsilon_2 R_{1,2} =\frac1N + \varepsilon_1\varepsilon_2R_{1,2}^-$, we can write
\begin{equation}\label{p4}
\nu((R_{1,2}^{-})^2)=\nu(R_{1,2}^2)-\frac2N\nu(\varepsilon_1\varepsilon_2R_{1,2})+\frac1{N^2}.
\end{equation}
Finally
\begin{eqnarray}
{\nu''\!}_0((R_{1,2}^{-})^2)&=&\beta^4\left[\nu_0((R_{1,2}^{-})^4)-
4\nu_0((R_{1,2}^{-})^2(R_{1,3}^{-})^2)+3\nu_0((R_{1,2}^{-})^2(R_{3,4}^{-})^2)
\right]\nonumber\\
&=&\beta^4\left[\nu(R_{1,2}^4)-4\nu(R_{1,2}^2R_{1,3}^2)+3\nu(R_{1,2}^2R_{3,4}^2)\right]+O\left(\frac1{N^{\frac52}}\right),\label{p5}
\end{eqnarray}
where the first equality follows from the computation of the derivative and the second one is again an easy consequence of the relation $\varepsilon_1 \varepsilon_2 R_{1,2} =\frac1N + \varepsilon_1\varepsilon_2R_{1,2}^-$.
Then, putting together (\ref{p3}), (\ref{p4}) and (\ref{p5}) and using $\nu(\varepsilon_1\varepsilon_2R_{1,2})=\nu(R_{1,2}^2)$, we obtain
\begin{eqnarray}
{\nu'\!}_0(\varepsilon_1\varepsilon_2R_{1,2}^{-})&=&\beta^2
\nu(R_{1,2}^2)-\frac{2\beta^2}N\nu(R_{1,2}^2)+\frac{\beta^2}{N^2}\nonumber\\&&-\frac12\beta^6\nu(R_{1,2}^4)+
2\beta^6\nu(R_{1,2}^2R_{1,3}^2)-\frac32\beta^6\nu(R_{1,2}^2R_{3,4}^2)+O\left(\frac1{N^{\frac52}}\right).\label{p6}
\end{eqnarray}

{\it Study of ${\nu''\!}_0(\varepsilon_1\varepsilon_2R_{1,2}^{-})$}. Using similar arguments, that is, computation of the derivatives, Taylor's expansions and the relation
$\varepsilon_1 \varepsilon_2 R_{1,2} =\frac1N + \varepsilon_1\varepsilon_2R_{1,2}^-$ we can easily prove that
\begin{eqnarray}
{\nu''\!}_0(\varepsilon_1\varepsilon_2R_{1,2}^{-})&=&-4\beta^4
\nu_0(R_{1,2}^{-}R_{1,3}^{-}R_{2,3}^{-})\nonumber\\&=&-4\beta^4\left[
\nu(R_{1,2}^{-}R_{1,3}^{-}R_{2,3}^{-})-{\nu'\!}_0(R_{1,2}^{-}R_{1,3}^{-}R_{2,3}^{-})\right]+O\left(\frac1{N^{\frac52}}\right)\nonumber\\
&=&-4\beta^4\left[\nu(R_{1,2}R_{1,3}R_{2,3})-\frac3N\nu(R_{1,2}R_{1,3}\varepsilon_2\varepsilon_3)\right]+O\left(\frac1{N^{\frac52}}\right)\nonumber\\
&=&-4\beta^4\nu(R_{1,2}R_{1,3}R_{2,3})+O\left(\frac1{N^{\frac52}}\right). \label{p7}
\end{eqnarray}
Notice only that we have used that
${\nu'\!}_0(R_{1,2}^{-}R_{1,3}^{-}R_{2,3}^{-})=0$ by Proposition \ref{der1} and Lemma \ref{simp1} and that by the symmetry property and again by Lemma \ref{simp1}
$$\vert \nu(R_{1,2}R_{1,3}\varepsilon_2\varepsilon_3) \vert =
\vert \nu(R_{1,2}R_{1,3}R_{2,3}) \vert \le \frac{K}{N^{\frac32}}.$$

{\it Study of $\nu^{(3)}_0(\varepsilon_1\varepsilon_2R_{1,2}^{-})$}. By similar arguments we can now obtain
\begin{eqnarray}
\nu_0^{(3)}(\varepsilon_1\varepsilon_2R_{1,2}^{-})&=&\beta^6\left[\nu_0((R_{1,2}^{-})^4)-12\nu_0((R_{1,2}^{-})^2(R_{2,3}^{-})^2)\right.\nonumber\\
&&\left.+9\nu_0((R_{1,2}^{-})^2(R_{3,4}^{-})^2)+36\nu_0(R_{1,2}^{-}R_{1,3}^{-}R_{2,4}^{-}R_{3,4}^{-})\right]\nonumber\\
&=&\beta^6\left[\nu(R_{1,2}^4)-12\nu(R_{1,2}^2R_{2,3}^2)
+9\nu(R_{1,2}^2R_{3,4}^2)\right.\nonumber\\
&&\left.+36\nu(R_{1,2}R_{1,3}R_{2,4}R_{3,4})\right]+O\left(\frac1{N^{\frac52}}\right).\label{p8}
\end{eqnarray}

Using the same ideas, it is also easy to check that
$\vert \nu(R_{1,2} R_{1,3} R_{2,4} R_{3,4} ) \vert \le \frac{K}{N^{\frac52}}$.
Then, putting together (\ref{p1}), (\ref{p2}), (\ref{p6}), (\ref{p7}) and (\ref{p8}), and
observing that some terms cancels out, we get
\begin{eqnarray*}
\nu(R_{1,2}^2)&=&\frac1N+\beta^2\nu(R_{1,2}^2)+\frac{\beta^2}{N^2}-\frac{2\beta^2}N\nu(\varepsilon_1\varepsilon_2R_{1,2})\\
&&-2\beta^4\nu(R_{1,2}R_{1,3}R_{2,3})-\frac13\beta^6\nu(R_{1,2}^4)+O\left(\frac1{N^{\frac52}}\right),
\end{eqnarray*}
or the equivalent equality
\begin{eqnarray}
\nu(R_{1,2}^2)&=&\frac{1}{1-\beta^2} \Big( \frac1N+ \frac{\beta^2}{N^2}-\frac{2\beta^2}N\nu(R^2_{1,2})\nonumber\\
&&-2\beta^4\nu(R_{1,2}R_{1,3}R_{2,3})-\frac13\beta^6\nu(R_{1,2}^4) \Big)+ O\left(\frac1{N^{\frac52}}\right).\label{expe1}
\end{eqnarray}

{\it Study of $\nu(R_{1,2}^4)$}. Applying the same type of arguments it is easy to get
\begin{eqnarray*}
\nu(R_{1,2}^4)&=&\nu(\varepsilon_1 \varepsilon_2 R_{1,2}^3)\\
&=&\nu(\varepsilon_1 \varepsilon_2 (R_{1,2}^-)^3)+ \frac{3}{N}\nu((R_{1,2}^-)^2)+ O\left(\frac1{N^{\frac52}}\right)\\
&=&\nu_0(\varepsilon_1 \varepsilon_2 (R_{1,2}^-)^3)+
{\nu'\!}_0(\varepsilon_1 \varepsilon_2 (R_{1,2}^-)^3)+
\frac{3}{N}\nu((R_{1,2}^-)^2)+ O\left(\frac1{N^{\frac52}}\right)\\
&=&\beta^2 \nu_0( R_{1,2}^4)+
\frac{3}{N}\nu((R_{1,2}^2))+ O\left(\frac1{N^{\frac52}}\right)\\
&=&\beta^2 \nu_0( R_{1,2}^4)+
\frac{3}{N^2(1-\beta^2)}+ O\left(\frac1{N^{\frac52}}\right),
\end{eqnarray*}
where we have used that (\ref{expe1}) implies that $\nu(R_{1,2}^2)=
\frac{1}{N(1-\beta^2)}+ O\left(\frac1{N^{\frac32}}\right)$. So it is clear that
\begin{equation}\label{exp2}
\nu(R_{1,2}^4)= \frac{3}{N^2(1-\beta^2)^2}+ O\left(\frac1{N^{\frac52}}\right),
\end{equation}

{\it Study of $\nu(R_{1,2}R_{1,3}R_{2,3})$}. Using the same ideas it is not difficult to obtain
\begin{equation}\label{exp3}
\nu(R_{1,2}R_{1,3}R_{2,3})= \frac{1}{N^2(1-\beta^2)^3}+ O\left(\frac1{N^{\frac52}}\right),
\end{equation}

Then, plugging (\ref{exp2}) and (\ref{exp3}) into (\ref{expe1}), and applying iteratively this new equality we get
\begin{eqnarray*}
\nu(R_{1,2}^2)&=&\frac1{N(1-\beta^2)}+\frac{\beta^2}{N^2(1-\beta^2)}-\frac{2\beta^2}N\frac1{N(1-\beta^2)^2}\\
&&-\frac{2\beta^4}{N^2(1-\beta^2)^4}-\frac{\beta^6}3\frac3{N^2(1-\beta^2)^3}+ O\left(\frac1{N^{\frac52}}\right)
\\
&=&\frac1{N(1-\beta^2)}-\frac{\beta^2+\beta^4}{N^2(1-\beta^2)^4}+ O\left(\frac1{N^{\frac52}}\right).
\end{eqnarray*}
\hfill $\Box$

\end{dem}

\section{Proof of the main result}

In this section we give the proof of Proposition \ref{propgen}. This proof is based on a
sequence of transformations that generalizes the tools used in the proof of the expansion of $\nu (R_{1,2}^2).$ We begin by introducing a new definition and explaining the four transformations we will use.

\begin{defin} We will say that a term is of order $m$ if its absolute value can be bounded by a constant times $\frac{1}{N^m}$.
\end{defin}

\begin{obs}\label{obsa} From Lemma \ref{simp1}, notice that the terms
 ${\displaystyle
\frac1{N^j}\nu(\prod_{i=1}^k\varepsilon_{j_i}\varepsilon_{l_i}R_{j_i,l_i})}$, ${\displaystyle
\frac1{N^j}\nu(\prod_{i=1}^k\varepsilon_{j_i}\varepsilon_{l_i}R_{j_i,l_i}^{-})}$
and ${\displaystyle
\frac1{N^j}\nu_0(\prod_{i=1}^k\varepsilon_{j_i}\varepsilon_{l_i}R_{j_i,l_i})}$
are of order $j+\frac k2$.
\end{obs}

Let us now introduce the four transformations that we will use througout this paper.

\bigskip

\noindent {\bf Transformation a.} {\it From $R$ to $R^-$.}

Using the relation
$\varepsilon_{j_i}\varepsilon_{l_i}R_{j_i,l_i}=\varepsilon_{j_i}\varepsilon_{l_i}R_{j_i,l_i}^{-}+\frac1N$,
we get
$$\nu(\prod_{i=1}^k\varepsilon_{j_i}\varepsilon_{l_i}R_{j_i,l_i})=\frac1{N^{k}}+
\nu(\prod_{i=1}^k\varepsilon_{j_i}\varepsilon_{l_i}R_{j_i,l_i}^{-})+T_a,$$
where $T_a$ is a finite sum of terms of the type
$$\frac1{N^m} \nu(\prod_{i=1}^{k-m}\varepsilon_{j_i}\varepsilon_{l_i}R_{j_i,l_i}^{-})$$
with $k-1 \ge m \ge 1$, that are of order greater than or equal to $\frac{k+1}2$ (see Remark \ref{obsa}).

\bigskip

\noindent {\bf Transformation b.} {\it Taylor's expansion.}

For a fixed $n \ge \frac{k}{2}$, using Taylor's formula, we obtain
$$\nu(\prod_{i=1}^k\varepsilon_{j_i}\varepsilon_{l_i}R_{j_i,l_i}^{-})=
\nu_0(\prod_{i=1}^k\varepsilon_{j_i}\varepsilon_{l_i}R_{j_i,l_i}^{-})
+\sum_{m=1}^{2n-k} \frac1{m!}
\nu_0^{(m)}(\prod_{i=1}^k\varepsilon_{j_i}\varepsilon_{l_i}R_{j_i,l_i}^{-})+O\left(\frac1{N^{n+\frac{1}{2}}}\right).$$
Applying Proposition \ref{deriv} $m$ times we get that
$\nu_0^{(m)}(\prod_{i=1}^k\varepsilon_{j_i}\varepsilon_{l_i}R_{j_i,l_i}^{-})$ is a finite sum
of terms of the form
\begin{equation}\label{eq1}
 \beta^{2m} K(m,k,I')
\nu_0(\prod_{i=1}^k\varepsilon_{j_i}\varepsilon_{l_i}R_{j_i,l_i}^{-}
\prod_{i'=1}^m\varepsilon_{j_{i'}}\varepsilon_{l_{i'}}R_{j_{i'}l_{i'}}^{-}),
\end{equation}
where $I'$ designates the sequence $\{(j_i,l_i), 1 \le i \le m \}$.

So we can write
$$\nu(\prod_{i=1}^k\varepsilon_{j_i}\varepsilon_{l_i}R_{j_i,l_i}^{-})=
\nu_0(\prod_{i=1}^k\varepsilon_{j_i}\varepsilon_{l_i}R_{j_i,l_i}^{-})+
{\nu'\!}_0(\prod_{i=1}^k\varepsilon_{j_i}\varepsilon_{l_i}R_{j_i,l_i}^{-})+
T_b+O\left(\frac1{N^{n+\frac{1}{2}}}\right),$$
where $T_b$ is a finite sum of terms of the type
$$K(\beta) \nu_0(\prod_{i=1}^{p}\varepsilon_{j_i}\varepsilon_{l_i}R_{j_i,l_i}^{-})$$
with $2n \ge p \ge k+2$. These terms are all of order greater than or equal to $\frac{k}{2} + 1$ (see Remark \ref{obsa}).

\bigskip

\noindent {\bf Transformation c.} {\it Inverse Taylor's expansion.}

Fixing again $n  \ge \frac{k}{2}$, using a Taylor's formula we get
$$\nu_0(\prod_{i=1}^k\varepsilon_{j_i}\varepsilon_{l_i}R_{j_i,l_i}^{-})=
\nu(\prod_{i=1}^k\varepsilon_{j_i}\varepsilon_{l_i}R_{j_i,l_i}^{-})
-\sum_{m=1}^{2n-k} \frac1{m!}
\nu_0^{(m)}(\prod_{i=1}^k\varepsilon_{j_i}\varepsilon_{l_i}R_{j_i,l_i}^{-})+O\left(\frac1{N^{n+\frac{1}{2}}}\right).$$
Using again the expression (\ref{eq1}) we can write
\begin{eqnarray*}
\nu_0(\prod_{i=1}^k\varepsilon_{j_i}\varepsilon_{l_i}R_{j_i,l_i}^{-})&=&
\nu(\prod_{i=1}^k\varepsilon_{j_i}\varepsilon_{l_i}R_{j_i,l_i}^{-})\\
&-&\sum_{m=1}^{2n-k} \sum_{I'} \beta^{2m} K(m,k)
\nu_0(\prod_{i=1}^k\varepsilon_{j_i}\varepsilon_{l_i}R_{j_i,l_i}^{-}
\prod_{j=1}^m\varepsilon_{j_{i'}}\varepsilon_{l_{i'}}R_{j_{i'}l_{i'}}^{-})+O\left(\frac1{N^{n+\frac{1}{2}}}\right).
\end{eqnarray*}
Repeating this procedure with all the terms
${\displaystyle
\nu_0(\prod_{i=1}^k\varepsilon_{j_i}\varepsilon_{l_i}R_{j_i,l_i}^{-}
\prod_{j=1}^m\varepsilon_{j_{i'}}\varepsilon_{l_{i'}}R_{j_{i'}l_{i'}}^{-})}$,
we obtain that
$$\nu_0(\prod_{i=1}^k\varepsilon_{j_i}\varepsilon_{l_i}R_{j_i,l_i}^{-})=
\nu(\prod_{i=1}^k\varepsilon_{j_i}\varepsilon_{l_i}R_{j_i,l_i}^{-})
+T_c+O\left(\frac1{N^{n+\frac{1}{2}}}\right),$$
where $T_c$ is a finite sum of terms of the type
$$K(\beta) \nu(\prod_{i=1}^{p}\varepsilon_{j_i}\varepsilon_{l_i}R_{j_i,l_i}^{-})$$
with $2n \ge p \ge k+1$ that are of  order greater than or equal to $\frac{k+1}{2}$ (see Remark \ref{obsa}).

\bigskip

\noindent {\bf Transformation d.} {\it From $R^-$ to $R$ }

Using the relation
$\varepsilon_{j_i}\varepsilon_{l_i}R_{j_i,l_i}^{-}=\varepsilon_{j_i}\varepsilon_{l_i}R_{j_i,l_i}-\frac1N$,
we get now that
$$\nu(\prod_{i=1}^k\varepsilon_{j_i}\varepsilon_{l_i}R_{j_i,l_i}^-)=\big(- \frac1{N} )^k+
\nu(\prod_{i=1}^k\varepsilon_{j_i}\varepsilon_{l_i}R_{j_i,l_i})+T_d,$$
where $T_d$ is a finite sum of terms of the type
$$\big( - \frac1{N} \big)^m \nu(\prod_{i=1}^{k-m}\varepsilon_{j_i}\varepsilon_{l_i}R_{j_i,l_i}^{-})$$
with $k-1 \ge m \ge 1$, that are of order greater than or equal to $\frac{k+1}2$ (see Remark \ref{obsa}).

\bigskip

The following remarks explain how can we apply these transformations.

\begin{obs}\label{tr2} For any fixed $n, j, k$ integers
applying {\sl transformation c} and {\sl transformation d} we can develop the
expression
$$
{\displaystyle
\frac1{N^j}\nu_0(\prod_{i=1}^k\varepsilon_{j_i}\varepsilon_{l_i}R_{j_i,l_i}^{-})}
$$
into a remainder of order $n + \frac{1}{2}$ plus a finite sum of terms of the type
$$K(\beta)\frac1{N^l}\nu(\prod_{i=1}^m\varepsilon_{j_i}\varepsilon_{l_i}R_{j_i,l_i})$$
with $n \ge l+\frac m2\geq j+\frac k2$, $l+m \ge j+k$, $l\geq j$. Notice again that all this new terms are of order greater than the initial one, except for the term
$$
{\displaystyle
\frac1{N^j}\nu(\prod_{i=1}^k\varepsilon_{j_i}\varepsilon_{l_i}R_{j_i,l_i})}
$$
that is of the same order.
\end{obs}

\begin{obs}
\label{tr3} Fixed $n, j, k$, using {\sl transformation b} and Remark \ref{tr2} we can expand the expression
$${\displaystyle
\frac1{N^j}\nu(\prod_{i=1}^k\varepsilon_{j_i}\varepsilon_{l_i}R_{j_i,l_i}^-)}
$$
into a remainder of order $n + \frac{1}{2}$ plus a finite sum of terms of the type
$$K(\beta)\frac1{N^l}\nu(\prod_{i=1}^m\varepsilon_{j_i}\varepsilon_{l_i}R_{j_i,l_i})$$
amb $n \ge l+\frac m2\geq j+\frac k2$, $l+m \ge j+k$, $l\geq j$. Notice also that all this new terms are of order greater than the initial one, except the term
${\displaystyle
\frac{K_0(\beta)}{N^j}\nu(\prod_{i=1}^k\varepsilon_{j_i}\varepsilon_{l_i}R_{j_i,l_i})}$ that is of the initial order. Moreover, if we denote by  $\lambda$ the $\varepsilon$-order of our initial expression
${\displaystyle
\frac1{N^j}\nu(\prod_{i=1}^k\varepsilon_{j_i}\varepsilon_{l_i}R_{j_i,l_i}^-)}$, we can check that $K_0(\beta)=\beta^2$ when $\lambda=0$ and $K_0(\beta)=0$ when $\lambda>0$.
\end{obs}

\begin{dem}{\bf [Proposition \ref{propgen}]}
Set $\lambda$ for the $\varepsilon$-order associated to
${\displaystyle
\prod_{i=1}^k\varepsilon_{j_i}\varepsilon_{l_i}R_{j_i,l_i}}.$
We need only to consider 3 cases.

\begin{itemize}
\item {\it Case a):} $\lambda\geq4$.

Applying {\sl Transformation a} we obtain that
$$\nu(\prod_{i=1}^k\varepsilon_{j_i}\varepsilon_{l_i}R_{j_i,l_i})=
\nu(\prod_{i=1}^k\varepsilon_{j_i}\varepsilon_{l_i}R_{j_i,l_i}^{-})+T_a^*,$$
where $T_a^*$ is a finite sum of terms of the type
$$\frac1{N^j}\nu(\prod_{i=1}^r\varepsilon_{j_i}\varepsilon_{l_i}R_{j_i,l_i}^{-})$$
with $j+\frac r2>\frac k2$. Notice that all these terms can be studied by Remark \ref{tr3}.

So, it suffices to deal with the term ${\displaystyle
\nu(\prod_{i=1}^k\varepsilon_{j_i}\varepsilon_{l_i}R_{j_i,l_i}^{-})}$.
Since $\lambda\geq4$ it is easily seen that
$$
\nu_0(\prod_{i=1}^k\varepsilon_{j_i}\varepsilon_{l_i}R_{j_i,l_i}^{-})=
{\nu'\!}_0(\prod_{i=1}^k\varepsilon_{j_i}\varepsilon_{l_i}R_{j_i,l_i}^{-})=0$$
and moreover, applying {\sl Transformation b} we get a finite sum of terms of the type
$$K(\beta)\nu_0(\prod_{i=1}^r\varepsilon_{j_i}\varepsilon_{l_i}R_{j_i,l_i}^{-})$$
with $r\geq k+2$. Remark \ref{tr2} completes the proof of this case.

\item {\it Case b):} $\lambda=2$.

Applying {\sl Transformation a} we can write
\begin{equation}\label{lam1}
\nu(\prod_{i=1}^k\varepsilon_{j_i}\varepsilon_{l_i}R_{j_i,l_i})=
\nu(\prod_{i=1}^k\varepsilon_{j_i}\varepsilon_{l_i}R_{j_i,l_i}^{-})+T_a^{**}+S_a^{**},
\end{equation}
where $T_a^{**}$ is a finite sum of terms of the type $$\frac1{N}
\nu(\prod_{i=1}^{k-1}\varepsilon_{j_i}\varepsilon_{l_i}R_{j_i,l_i}^{-})$$
and $S_a^{**}$ is a finite sum of terms of the type
$$\frac1{N^j}\nu(\prod_{i=1}^r\varepsilon_{j_i}\varepsilon_{l_i}R_{j_i,l_i}^{-})$$
with $j>1$, $r+j=k$, $\frac r2+j\geq\frac k2+1$.

Let us study first the term ${\displaystyle
\nu(\prod_{i=1}^k\varepsilon_{j_i}\varepsilon_{l_i}R_{j_i,l_i}^{-})}$.
Since $\lambda=2$ a trivial verification shows that ${\displaystyle
\nu_0(\prod_{i=1}^k\varepsilon_{j_i}\varepsilon_{l_i}R_{j_i,l_i}^{-})}=0$
and that there exists $j_{k+1},l_{k+1}$ verifying
$$
{\nu'\!}_0(\prod_{i=1}^k\varepsilon_{j_i}\varepsilon_{l_i}R^-_{j_i,l_i})=
\beta^2
\nu_0((\prod_{i=1}^k\varepsilon_{j_i}\varepsilon_{l_i}R_{j_i,l_i}^{-})(\varepsilon_{j_{k+1}}\varepsilon_{l_{k+1}}R_{j_{k+1},l_{k+1}}^{-})),$$
and such that the $\varepsilon$-order of
${\displaystyle(\prod_{i=1}^k\varepsilon_{j_i}\varepsilon_{l_i}R_{j_i,l_i}^{-})(\varepsilon_{j_{k+1}}\varepsilon_{l_{k+1}}R_{j_{k+1},l_{k+1}}^{-})}$ is zero.
So, applying {\sl Transformation b} we can write
\begin{equation}\label{lam2}
\nu(\prod_{i=1}^k \varepsilon_{j_i} \varepsilon_{l_i} R_{j_i,l_i}^{-})
= \beta^2
\nu_0((\prod_{i=1}^k\varepsilon_{j_i}\varepsilon_{l_i}R_{j_i,l_i}^{-
})(\varepsilon_{j_{k+1}}\varepsilon_{l_{k+1}}R_{j_{k+1},l_{k+1}}^{-}))+T_b^{***}+O\left(\frac1{N^{n+\frac{1}{2}}}\right),
\end{equation}
where $T_b^{***}$  is a finite sum of terms of the type
$$K(\beta)
\nu_0(\prod_{i=1}^r\varepsilon_{j_i}\varepsilon_{l_i}R_{j_i,l_i})$$
with $ r \geq k+2$.

Let us study now the term ${\beta^2
\nu_0((\prod_{i=1}^k\varepsilon_{j_i}\varepsilon_{l_i}R_{j_i,l_i}^{-
})(\varepsilon_{j_{k+1}}\varepsilon_{l_{k+1}}R_{j_{k+1},l_{k+1}}^{-}))}$.
Notice first that
 $
\nu_0'((\prod_{i=1}^k\varepsilon_{j_i}\varepsilon_{l_i}R_{j_i,l_i}^{-
})(\varepsilon_{j_{k+1}}\varepsilon_{l_{k+1}}R_{j_{k+1},l_{k+1}}^{-}))=0$. Then
 applying {\sl Transformations c} and {\sl d} we can write
\begin{multline}
\beta^2
\nu_0((\prod_{i=1}^k\varepsilon_{j_i}\varepsilon_{l_i}R_{j_i,l_i}^{-
})(\varepsilon_{j_{k+1}}\varepsilon_{l_{k+1}}R_{j_{k+1},l_{k+1}}^{-}))
\\
= \beta^2
\nu((\prod_{i=1}^k\varepsilon_{j_i}\varepsilon_{l_i}R_{j_i,l_i})(\varepsilon_{j_{k+1}}\varepsilon_{l_{k+1}}R_{j_{k+1},l_{k+1}}))+T_{cd}^{*}+O\left(\frac1{N^{n+\frac{1}{2}}}\right),
\label{lam3}\end{multline}
where $T_{cd}^{*}$  is a sum finite of terms of the type
$$K(\beta)\frac1{N^j}
\nu_0(\prod_{i=1}^r\varepsilon_{j_i}\varepsilon_{l_i}R_{j_i,l_i})$$
with $n \ge j+\frac r2>\frac{k+2}2$.

Finally, by a symmetry argument it is clear that
\begin{eqnarray}
\nu((\prod_{i=1}^k\varepsilon_{j_i}\varepsilon_{l_i}R_{j_i,l_i})(\varepsilon_{j_{k+1}}\varepsilon_{l_{k+1}}R_{j_{k+1},l_{k+1}}))=
\nu((\prod_{i=1}^k R_{j_i,l_i})R_{j_{k+1},l_{k+1}})\nonumber \\
= \nu(\varepsilon_{j_{k+1}}\varepsilon_{l_{k+1}}\prod_{i=1}^k R_{j_i,l_i})
= \nu(\prod_{i=1}^k\varepsilon_{j_i}\varepsilon_{l_i}R_{j_i,l_i}).
\label{lam4}\end{eqnarray}

Putting together (\ref{lam1}), (\ref{lam2}), (\ref{lam3}) and (\ref{lam4}) we obtain
$$\nu(\prod_{i=1}^k\varepsilon_{j_i}\varepsilon_{l_i}R_{j_i,l_i})=\beta^2\nu(\prod_{i=1}^k\varepsilon_{j_i}\varepsilon_{l_i}R_{j_i,l_i})+T_a^{**}+S_a^{**}+T_b^{***}+T_{cd}^{*}+O\left(\frac1{N^{n+\frac{1}{2}}}\right),$$
or equivalently
$$\nu(\prod_{i=1}^k\varepsilon_{j_i}\varepsilon_{l_i}R_{j_i,l_i})=\frac{1}{1-\beta^2} \Big( T_a^{**}+S_a^{**}+T_b^{***}+T_{cd}^{*} \Big) +O\left(\frac1{N^{n+\frac{1}{2}}}\right).$$

Applying Remark \ref{tr3} to the terms of $T_a^{**}$ and $S_a^{**}$, and Remark \ref{tr2} to the terms of $T_b^{***}$ it is clear that we only obtain a sum of terms of the type
$$K(\beta)\frac1{N^j}
\nu(\prod_{i=1}^r\varepsilon_{j_i}\varepsilon_{l_i}R_{j_i,l_i})$$
with $j+\frac r2\geq \frac{k+1}2$.

Finally observe that the terms satisfying $j+\frac r2 = \frac{k+1}2$ come from $T_a^{**}$ and so they have $j=1$.

\item {\it Case c):} $\lambda=0$.

By symmetry we have
$$\nu(\prod_{i=1}^k\varepsilon_{j_i}\varepsilon_{l_i}R_{j_i,l_i})=
\nu(\prod_{i=1}^{k-1}\varepsilon_{j_i}\varepsilon_{l_i}R_{j_i,l_i}).$$

The expression $\prod_{i=1}^{k-1}\varepsilon_{j_i}\varepsilon_{l_i}R_{j_i,l_i}$
is of $\varepsilon$-order equal to 2, so we can reduce this case to the case b).
Notice that since we have now order $\frac{k-1}{2}$, applying the case b), we will obtain terms of order $\frac{k}{2}$ , but with $j>0$.

\end{itemize}

\hfill $\Box$

\end{dem}

\section{Computer Program}

In the proof given in the above section we found the basis of the algorithm that allows us to compute the expansion of a term ${\displaystyle\nu(\prod_{i=1}^k
\varepsilon_{j_i}\varepsilon_{l_i}R_{j_i,l_i})}$.

The algorithm is the following

\begin{itemize}
\item {1.} Study of $\varepsilon$-order.  When the $\varepsilon$-order is equal to 0 or to 2 we keep away part of the term ${\displaystyle-\frac{\beta^2}{1-\beta^2}\nu(\prod_{i=1}^k
\varepsilon_{j_i}\varepsilon_{l_i}R_{j_i,l_i})}$ that will be simplified later.
When the $\varepsilon$-order is equal to 0 we also use a symmetry argument.

\item {2.} From $R_{j,l}$ to $R_{j,l}^-$. (transformation a).

\item {3.} From ${\displaystyle\nu(\prod_{i=1}^k
\varepsilon_{j_i}\varepsilon_{l_i}R_{j_i,l_i}^-)}$ to
${\displaystyle\nu_0(\prod_{i=1}^k
\varepsilon_{j_i}\varepsilon_{l_i}R_{j_i,l_i}^-)}$. (transformation b)

\item {4.} From ${\displaystyle\nu_0(\prod_{i=1}^k
\varepsilon_{j_i}\varepsilon_{l_i}R_{j_i,l_i}^-)}$ to
${\displaystyle\nu(\prod_{i=1}^k
\varepsilon_{j_i}\varepsilon_{l_i}R_{j_i,l_i}^-)}$. (transformation c)

\item {5.} From $R_{j,l}^-$ to $R_{j,l}$. (transformation d).
\end{itemize}

After each one of these steps we have a procedure that simplifies the terms obtained. That is, we eliminate the terms with greater order than the grade of our expansion.
In order to do it we give a refinement of Lemma \ref{simp1}
that allows us to simplify more terms.

\begin{prop}\label{simp2}
Given ${\displaystyle\prod_{i=1}^k
\varepsilon_{j_i}\varepsilon_{l_i}R_{j_i,l_i}}$, set $\lambda$ its corresponding $\varepsilon$-order. Then
\begin{eqnarray*}
&(a)&
|\nu(\prod_{i=1}^k\varepsilon_{j_i}\varepsilon_{l_i}R_{j_i,l_i}^{-})|\leq\frac{K}{N^{\frac
k2 +\frac \lambda4}},\\
&(b)&|\nu(\prod_{i=1}^k\varepsilon_{j_i}\varepsilon_{l_i}R_{j_i,l_i})|\leq\frac{K}{N^{\frac
k2 +\frac \lambda4}}.
\end{eqnarray*}
\end{prop}
\begin{dem}

The case $\lambda=0$ has been studied in Lemma \ref{simp1}. Assume now $\lambda>0$.

Let us study first (a). Using a Taylor's expansion we get
\begin{eqnarray*}
\nu(\prod_{i=1}^k\varepsilon_{j_i}\varepsilon_{l_i}R_{j_i,l_i}^{-})&=&\nu_0(\prod_{i=1}^k\varepsilon_{j_i}\varepsilon_{l_i}R_{j_i,l_i}^{-})\\
&&+\frac1{n!}\sum_{n=1}^{\frac \lambda2
-1}\nu_0^{(n)}(\prod_{i=1}^k\varepsilon_{j_i}\varepsilon_{l_i}R_{j_i,l_i}^{-})+
\frac1{\left(\frac \lambda2\right)!}\nu_t^{\left(\frac \lambda2
\right)}(\prod_{i=1}^k\varepsilon_{j_i}\varepsilon_{l_i}R_{j_i,l_i}^{-}),
\end{eqnarray*}
for some $t \in [0,1]$. Notice that all the terms except the last one have strictly positive $\varepsilon$-order. So,
from Lemma \ref{simp1} we easily get
$$\nu(\prod_{i=1}^k\varepsilon_{j_i}\varepsilon_{l_i}R_{j_i,l_i}^{-})=\frac1{\left(\frac \lambda2\right)!}\nu_t^{\left(\frac \lambda2
\right)}(\prod_{i=1}^k\varepsilon_{j_i}\varepsilon_{l_i}R_{j_i,l_i}^{-}).$$
From Proposition \ref{der1} and Lemma \ref{simp1} we get
\begin{eqnarray*}
|\nu_t^{\left(\frac \lambda2
\right)}(\prod_{i=1}^k\varepsilon_{j_i}\varepsilon_{l_i}R_{j_i,l_i}^{-})|&\leq&\frac{K}{N^{\frac
\lambda4}}\nu\left(\left(\prod_{i=1}^k\varepsilon_{j_i}\varepsilon_{l_i}R_{j_i,l_i}^{-}\right)^2\right)^{\frac12}\\
&\leq&\frac{K}{N^{\frac \lambda4}}\left(\frac{K}{N^{\frac
{2k}2}}\right)^{\frac12}\\&=&\frac{K}{N^{\frac \lambda4 +\frac k2}}.
\end{eqnarray*}

To deal with (b) we will use the relation
$\varepsilon_{j_i}\varepsilon_{l_i}R_{j_i,l_i}=\varepsilon_{j_i}\varepsilon_{l_i}R_{j_i,l_i}^{-}+\frac1N$.
Then
$$\nu(\prod_{i=1}^k\varepsilon_{j_i}\varepsilon_{l_i}R_{j_i,l_i}^{-})=\sum_{j=0}^k\left(\frac1{N^{j}}\right)
\sum_{I}
\nu(\prod_{i=1}^{k-j}\varepsilon_{j_i}\varepsilon_{l_i}R_{j_i,l_i}^{-}).$$
Notice that all the terms
${\displaystyle\prod_{i=1}^{k-j}\varepsilon_{j_i}\varepsilon_{l_i}R_{j_i,l_i}^{-}}$
are of $\varepsilon$-order greater than or equal to $\max(0,\lambda-2j)$.
So, using the inequality obtained for (a) we have
$$
\left|\frac1{N^j}\nu(\prod_{i=1}^{k-j}\varepsilon_{j_i}\varepsilon_{l_i}R_{j_i,l_i}^{-})\right| \leq \frac1{N^j}\frac{K}{N^{\frac{k-j}2}+\frac{\lambda^\alpha}4} \le
\frac{K}{N^{\frac k2 + \frac \lambda4}}.
$$
\hfill $\Box$

\end{dem}

\smallskip

In the Appendix you can find the main parts of a computer program where we have developped this algorithm. This
computer program has been written in Maple7. The program is of free use and is available from the authors.  One must introduce the expression
$\displaystyle{\nu(\prod_{i=1}^k\varepsilon_{j_i}\varepsilon_{l_i}R_{j_i,l_i})}$
using the notation $[[j_1,l_1],[j_2,l_2],\ldots,[j_k,l_k]]$. Then to call the computer program and to obtain the expansion of such expression until order m, one has to write
$tot([[[j_1,l_1],[j_2,l_2],\ldots,[j_k,l_k]]],m).$

\begin{obs}
For instance we can get
\begin{eqnarray*}
\nu(\varepsilon_2 \varepsilon_2 R_{1,2} R_{1,3}) &=& \frac{1}{(1-\beta^2)^3 N^2} +
O\left(\frac1{N^3}\right),\,\,\,\,\,tot([[[1, 2],[1, 3]]],2);\\
\nu( R_{1,2} R_{1,3}R_{2,3} R_{2,4}) &=& \frac{1+\beta}{(1-\beta^2)^5 N^3} +
O\left(\frac1{N^4}\right),\,\,\,\,\,tot([[[1, 2],[1, 3],[2,4],[3,4]]],3);\\
\nu(\varepsilon_1 \varepsilon_2 \varepsilon_3 \varepsilon_4 R_{1,2} R_{3,4}) &=& \frac{1}{(1-\beta^2)^2 N^2} +
O\left(\frac1{N^3}\right),\,\,\,\,\,tot([[[1, 2],[3, 4]]],2).
\end{eqnarray*}
\end{obs}

\section*{Appendix}

{\small {\small {\small {\small

\begin{itemize}



\item[]  {\bf Main program (tot)}
Parameters: {\it expres} indicates the expression we want to expand
and {\it derivord} is the order of the expansion

\item[] $>$ tot:=proc(expres,derivord)

\begin{itemize}

\item[]  local taula,coef,i,k,m;

\item[]  global resfin,resfin2,kdes,hnm;

\begin{itemize}

\item[]  taula := expres;
 coef := [1];
 resfin:=0;
resfin2:=0;
kdes:=derivord;
hnm:=1;

\item[]  maxima(taula[1]);
k:=0;

\item[]     for i from 1 to r do; 

\begin{itemize}

\item        m:=(numboccur(taula[j],i) mod 2);
if m=1 then  k:=k+1; end if;

\end{itemize}

\item[] end do;

\item[]  if kdes $<$ ((1/2)*nops(taula[1]))+k/2
 then hnm:=0; resfin2:=0;
end if;


\item[]  while (nops(taula)*hnm)$>$0 do

\begin{itemize}

\item[]     inici (taula,coef);
taula:=tauls;
coef:=coefs;
resfin:=collect((resfin+resu),1/N);

\end{itemize}

\item[]  end do;

\item[]  for i from 1 to kdes do
 resfin2:=resfin2+(simplify(coeff(resfin,N,-i))*N\^\,(-i));
 end do;

\item[]  print (RESULT, resfin2);

\end{itemize}

\item[]  end proc;

\end{itemize} 


\item[] {\bf Procedure inici}
 {\it This procedure calls the four transformations explained in the
paper}


\item[] {\bf Procedure rseparar}
{\it We study the value of the $\varepsilon$-order associated to each term. When
it is equal to 0 or to 2 we keep away a part of the term that
will be simplified later.}

\item[] {\bf Procedure maxima}
{\it Compute the maximum of a list}

\item[] {\bf Procedure rmesamenys} {\it Transformation a}

\item[] $>$ rmesamenys:=proc(taull,coefl)

\begin{itemize}

\item[]  local
i,j,ll,k,l,m,kter,long,coefll,coef,taula,taulab,fer,coefb;

\item[]  global taulabb,coefbb,resb;

\begin{itemize}

\item[]  taula:=taull;
coef:=coefl;
resb:=0;
taulab:=[$\,$];
coefb:=[$\,$];


\item[]  for i from 1 to nops(taula) do;

\begin{itemize}

\item[]      long:= nops(taula[i]);

\item[]      for j from 0 to 2\^\,long-1 do;

\begin{enumerate}

\item[]          kter:=0;
conversio(j,long);
ll:=[$\,$];
coefll:=coef[i];

\item[]          for k from 1 to long do;

\begin{enumerate}

\item[]              if jbase[k]=0 then coefll:=coefll*(1/N);

\item[]                            else
kter:=kter+(1/2);
ll:=[op(ll),taula[i][k]];
end if;

\end{enumerate}

\item[]          end do;

\item[]          kter:=kter-degree(coefll,N)+(comptar(ll)/4);

\item[]          if kter$<$kdes+(1/2) and ll$<>$[$\,$] then

\begin{enumerate}

\item[]                     taulab:=[op(taulab),ll];
coefb:=[op(coefb),coefll];

\end{enumerate}

\item[]          end if;

\item[]          if kter$<$kdes+(1/2) and ll=[$\,$] then
resb:=resb+coefll;
end if;

\end{enumerate}

\item[]       end do;

\end{itemize}

\item[]  end do;


\item[]  taulabb := [$\,$];
coefbb := [$\,$];
taulabb := [op(taulabb),taulab[1]];

\item[] coefbb := [op(coefbb),coefb[1]];


\item[]  for i from 2 to nops(taulab) do

\begin{itemize}

\item[]    for j from 1 to nops(taulabb) while
(taulab[i]$<>$taulabb[j] or
\item[] degree(coefb[i],N)$<>$degree(coefbb[j],N)) do

\begin{enumerate}

\item[]         if j=nops(taulabb) then taulabb :=
[op(taulabb), taulab[i]];

\begin{enumerate}

\item[]                                 coefbb :=
[op(coefbb), 0];

\end{enumerate}

\item[]         end if;

\end{enumerate}

\item[]    end do;

\item[]    fer:= coefbb[j]+coefb[i];
coefbb:=subsop(j=fer,coefbb);

\end{itemize}

\item[]  end do;

\end{itemize}

\item[]  end proc;

\end{itemize}


\item[] {\bf Procedure conversio} {\it Convert a number to a binary expression given in a list}

\item[] {\bf Procedure comptar} {\it It counts the ocurrences of each number in a list}

\item[] {\bf  Procedure simplif} {\it It joins all the terms with the same expression}

\item[] {\bf Procedure vdesenvola}
{\it Transformation b}

\item[] $>$ vdesenvola:=proc(taull,coefl)

\begin{itemize}

\item[]  local taula,coef,taula2,taula3,vtaula4,vcoef4,coef2,coef3,
\item[] vtaula3,vcoef2,vcoef3,coef5,coeft,ll,i,l,j,ik,k,x,kh,kter;

\item[]  global tauldh, coefdh;

\begin{itemize}

\item[]  tauldh:=[$\,$];
coefdh:=[$\,$];
kh:=kdes;
taula:=taull;
coef:=coefl;
taula2:=[$\,$];
coef2:=[$\,$];


\item[]  for i from 1 to nops(taula) do;

\begin{itemize}

\item[] kter:=2*(kh-((nops(taula[i])/2)-degree(coef[i],N)));

\item[]      vtaula3:=$<<0>>$;
 vcoef3:=$<<$coef[i]$>>$;

\item[]      for ik from 1 to nops(taula[i]) do;
vtaula3:=$<<$vtaula3 $\vert$ taula[i][ik]$>>$;
end do;

\item[] vtaula3:=SubMatrix(vtaula3,1..1,2..nops(taula[i])+1);

\item[]      vtaula4:=vtaula3;
vcoef4:=vcoef3;

\item[]      vanulareps(vtaula4,vcoef4,1);

\item[]      tauldh:=[op(tauldh),op(taulde)];
coefdh:=[op(coefdh),op(coefde)];

\item[]      for j from 1 to kter do;

\begin{enumerate}

\item[]          vderivar(vtaula3,vcoef3);

\item[]          vtaula3:=vtauld;
vcoef3:=vcoefd;

\item[]          vanulareps(vtauld,vcoefd,j);

\item[]          tauldh:=[op(tauldh),op(taulde)];
coefdh:=[op(coefdh),op(coefde)];

\end{enumerate}

\item[]      end do;

\end{itemize}

\item[]  end do;

\end{itemize}

\item[]  end proc;

\end{itemize}


\item[] {\bf Procedure vanulareps}
{\it It cancels all the terms that are null or that are of order
greater that our expansion (using Proposition \ref{simp2})}

\item[] {\bf Procedure vderivar}
{\it It gives us the derivative of an expression}

\item[] {\bf Procedure mmaxima} {\it Compute the maximum of a table}

\item[] {\bf Procedure vdesenvolb} {\it
Transformation c}

\item[] $>$ vdesenvolb:=proc(taull,coefl)

\begin{itemize}

\item[] local taula,taula2,taula3,taula4,tauldh,coefdh,ik,vtaula4,
\item[] vcoef4,coef,coef2,coef4,coef5,ll,i,l,j,k,x,kj,kter,n;

\item[]  global tauldj, coefdj;

\begin{itemize}

\item[]  tauldj:=[$\,$];
coefdj:=[$\,$];
taula4:=[$\,$];
taula2:=[$\,$];
coef2:=[$\,$];
coef4:=[$\,$];

\item[]  kj:=kdes;
taula:=taull;
coef:=coefl;


\item[]  for n from 1 while nops(taula)$>$0 do;

\begin{itemize}

\item[]     tauldh:=[$\,$];
coefdh:=[$\,$];

\item[]     tauldj:=[op(tauldj),op(taula)];
coefdj:=[op(coefdj),op(coef)];

\item[]     taula3:=taula;
taula:=[$\,$];

\item[]     for i from 1 to nops(taula3) do;

\begin{enumerate}

\item[] kter:=2*(kj-((nops(taula3[i])/2)-degree(coef[i],N)));

\item[]         vtaula4:=$<<0>>$;
vcoef4:=$<<$-coef[i]$>>$;

\item[]         for ik from 1 to nops(taula3[i]) do;
vtaula4:=$<<$vtaula4$\vert$taula3[i][ik]$>>$;
end do;

\item[] vtaula4:=SubMatrix(vtaula4,1..1,2..nops(taula3[i])+1);

\item[]         for j from 1 to kter do;

\begin{enumerate}

\item[]                vderivar(vtaula4,vcoef4);

\item[]                vtaula4:=vtauld;
vcoef4:=vcoefd;

\item[]                vanulareps(vtauld,vcoefd,j);

\item[]                tauldh:=[op(tauldh),op(taulde)];
coefdh:=[op(coefdh),op(coefde)];

\end{enumerate}

\item[]         end do;

\item[]         taula4:=[$\,$];
coef4:=[$\,$];

\end{enumerate}

\item[]     end do;

\item[]     taula:=tauldh;
coef:=coefdh;

\end{itemize}

\item[]  end do;

\end{itemize}

\item[]  end proc;

\end{itemize}


\item[] {\bf Procedure rmenysames} {\it Transformation d}

\end{itemize}

}}}}

\end{document}